\documentclass[11pt]{article}
\usepackage{amsfonts,latexsym,amsmath,amscd,geometry}
\usepackage{amssymb}
\usepackage{latexsym}
\usepackage{xcolor}

\newcommand \nc{\newcommand}
\newtheorem{theorem}{Theorem}[section]
\newtheorem{lemma}[theorem]{Lemma}

\newtheorem{corollary}[theorem]{Corollary}

\nc{\ba}{\begin{array}}\nc{\ea}{\end{array}}
\nc{\be}{\begin{eqnarray}}\nc{\ee}{\end{eqnarray}}
\nc{\beq}{\begin{equation}}\nc{\eeq}{\end{equation}}
\nc{\bex}{\begin{eqnarray*}}\nc{\eex}{\end{eqnarray*}}
\nc{\btm}{\begin{theorem}} \nc{\etm}{\end{theorem}}
\nc{\blm}{\begin{lemma}} \nc{\elm}{\end{lemma}}
\nc{\R}{\mathbb{R}}
\nc{\va}{\varphi}
\nc{\ve}{\varepsilon}
\def\S{\mathbb{S}}
\def\N{\mathcal{N}}
\def\z{\bar{z}}
\def\pf{\noindent{\bf Proof.\quad}}
\def\endpf{\hfill$\Box$}

\def\tr{\mbox{tr\,}}
\def\di{\mbox{div\,}}

\begin{document}
\title{Long-time dynamics of Ericksen-Leslie system on $\mathbb S^2$}

\author{Tao Huang\footnote{Department of Mathematics, Wayne State University, Detroit, MI 48202, USA.}  \quad Chengyuan Qu\thanks{Corresponding
author. E-mail: mathqcy@163.com} \footnote{School of Science, Dalian Minzu University, Dalian, Liaoning 116600, PRC.}\\
}
\date{}
\maketitle

\begin{abstract}
In this paper, we study the long-time behavior of full Ericksen-Leslie system modeling the hydrodynamics of nematic liquid crystals between two dimensional unit spheres. Under a weaker assumption for Leslie's coefficients, we give the key energy inequality for the global weak solution. At last,  inspired by the conditions on the simplified system, we establish several sufficient conditions which guarantee the uniform convergence of the system in $L^2$ and $H^k$ spaces as time tends to infinity under small initial data.

\end{abstract}

\section{Introduction}

Nematic liquid crystals are composed of rod-like molecules characterized by average alignment of  the long axes of neighboring molecules, which have simplest structures among various types of liquid crystals. The dynamic theory of nematic liquid crystals has been first proposed by Ericksen \cite{ericksen62} and Leslie \cite{leslie68} in the 1960's, which is a macroscopic continuum description of the time evolution of both flow velocity field and orientation order parameter of rod-like liquid crystals.

In this paper, we consider the following Ericksen-Leslie system on $\S^2\times(0,\infty)$, where $\S^2\subset\R^3$ is the unit sphere with standard metric $g_0$,
\begin{equation}\label{EL}
\begin{cases}
u_t+u\cdot\nabla u+\nabla P=-\di\big(\nabla d\odot\nabla d\big)+\di (\sigma^L(u,d)),\\
\di u=0, \\
\lambda_1(d_t+u\cdot\nabla d-\Omega d)+\lambda_2Ad=\Delta d+|\nabla d|^2d+\lambda_2(d^TAd)d,
\end{cases}
\end{equation}
where
$u(x,t): \S^2\times(0,\infty)\rightarrow T_x\S^2$ is the fluid velocity field, $d(x,t):\S^2\times(0,\infty)\rightarrow \S^2$ is the orientation order parameters of nematic material, $P(x,t): \S^2\times(0,\infty)\rightarrow \R$ is the pressure function.
For simplicity, we use $\nabla$, $\di$ and $\Delta$ to denote the gradient operator, divergence operator and the Laplace Beltrami operator on $\S^2$ with standard metric $g_0$, respectively. The convection term $u\cdot\nabla u$ in the first equation is the directed differentiation of $u$ with respect to the direction $u$ itself, which is interpreted as the covariant derivative $D_uu$. Here $D$ denotes the covariant derivative operator on $\S^2$.
Denote $d=(d_1, d_2, d_3)\in\S^2\hookrightarrow \R^3$ and $u=(u_1, u_2, u_2)\in T\S^2\hookrightarrow \R^3$. Then
$$
\big(\nabla d\odot\nabla d\big)_{ij}
=\langle\frac{\partial d}{\partial x_i}, \frac{\partial d}{\partial x_j}\rangle, \quad (u\cdot \nabla d)_i=g_0(u,\nabla d_i),\quad |\nabla d|^2=g_0(\nabla d, \nabla d),
$$
$\langle\cdot,\cdot\rangle$ denotes the inner product in $L^2$. Denote
$$
A_{ij}=\frac{1}{2}\left(\frac{\partial u_j}{\partial x_i}+\frac{\partial u_i}{\partial x_j}\right),\quad \Omega_{ij}=\frac{1}{2}\left(\frac{\partial u_i}{\partial x_j}-\frac{\partial u_j}{\partial x_i}\right),\quad N=d_t+u\cdot \nabla d-\Omega d,
$$
as the rate of strain tensor, skew-symmetric part of the strain rate, and the
rigid rotation part of direction changing rate by fluid vorticity, respectively.  The left side of third equation in \eqref{EL} is the kinematic transport, which represents the effect of the macroscopic flow field on the microscopic structure. The material coefficients $\lambda_1$ and $\lambda_2$ reflect the molecular shape and the slippery part between fluid and particles. The term with $\lambda_1$ represents the rigid rotation of molecules, while the term with $\lambda_2$ stands for the stretching of molecules by the flow. The viscous (Leslie) stress tensor $\sigma$ has the following form (cf. \cite{Les79})
$$
\sigma_{ij}^L(u,d)=\mu_1 d_kd_pA_{kp}d_id_j+\mu_2 N_id_j+\mu_3d_iN_j+\mu_4A_{ij}+\mu_5A_{ik}d_kd_j+\mu_6 A_{jk}d_kd_i
$$
The viscous coefficients $\mu_i$, $i=1,\ldots, 6$ are called Leslie's coefficients. The following relations are often assumed in the literature
\beq\label{parodi}
\mu_2+\mu_3=\mu_6-\mu_5
\eeq
\beq\label{nec1}
\lambda_1=\mu_3-\mu_2>0, \quad \lambda_2=\mu_6-\mu_5
\eeq
\beq\label{mus}
\mu_4>0,\quad 2\mu_1+3\mu_4+2\mu_5+2\mu_6>0,
\quad 2\mu_4+\mu_5+\mu_6>\frac{\lambda_2^2}{\lambda_1}
\eeq
The first relation is called Parodi's relation, which has been derived from Onsager reciprocal relations expressing the equality of certain relations between flows and forces in thermodynamic systems out of equilibrium (cf. \cite{Parodi70}).
The second two relations are the compatibility conditions. The third empirical relations are necessary to obtain the energy inequality (cf. \cite{Les79}).


In dimension two, the existence of global weak solutions to the initial and boundary value problem of \eqref{EL} has been constructed in \cite{huanglinwang14} for any bounded smooth domain. The weak solution has been proved to be regular except for finitely many times (see also \cite{hongxin12} for related results). The uniqueness of such weak solution has been proved in \cite{LTX16} and \cite{ wangwang14}. In dimension three, the global well-posedness combining with long time behaviors for the system \eqref{EL} around equilibria under various assumptions on the Leslie coefficients has been studied in \cite{ hieberpruss17, WZZ13, WXL13}.

There is a simplified system that has been first proposed in \cite{lin89} by neglecting the Leslie stress. There have been many works on the existence and partial regularity of the simplified system (see e.g. \cite{HLLW16}, \cite{llwwz19}, \cite{linlinwang10}, \cite{linwang16}, \cite{linwang10}). For the simplified system on the unit sphere $\S^2$, the uniform convergence of the solution in $L^2$ to a steady state exponentially as $t$ tends to infinity has been proved in \cite{wangxu14}.

In order to study the system \eqref{EL} on the unit sphere, we need to consider the initial data
\beq\label{initial}
u(x,0)=u_0,\quad d(x,0)=d_0
\eeq
with $u_0\in L^2(\S^2,T\S^2)$, $d_0\in H^1(\S^2, \S^2)$ and
\beq\label{u0con}
\int_{\S^2}u_0(x)\,dv_{0}=0,\quad \di u_0=0.
\eeq

Under the weaker  parameter relationship, the first main result concerns the existence and partial regularity of the global weak solutions to the system \eqref{EL} as follows.

\begin{theorem}\label{thm1}
Let any $u_0\in L^2(\S^2,T\S^2)$, $d_0\in H^1(\S^2, \S^2)$ satisfying \eqref{u0con}, then there exists a unique global weak solution $u\in L^{\infty}([0,\infty), L^2(\S^2,T\S^2))\cap L^{2}([0,\infty), H^1(\S^2,T\S^2))$, $d\in L^{\infty}([0,\infty), H^1(\S^2,\S^2))$ to the initial value problem \eqref{EL} and \eqref{initial} that satisfies the following properties
\begin{itemize}

\item[(1)] There exists an integer $L>0$ depending only on $u_0$ and $d_0$ such that $u, d\in C^\infty\left(\S^2\times\left((0,\infty)\setminus\{T_l\}_{l=1}^L\right)\right)$ for $0<T_1<\cdots<T_L<\infty$. At each singular time $T_l$, it holds for any $r>0$
\beq\label{blowupcr}
\limsup\limits_{t\uparrow T_l}\max\limits_{x\in\S^2}\int_{\S^2\cap B_r(x)}\left(|u|^2+|\nabla d|^2\right)(y,t)\,dv_0(g)(y)\geq 8\pi.
\eeq

\item[(2)] There exist a time sequence $t_i\rightarrow \infty$, a harmonic map $d_{\infty}\in C^\infty(\S^2,\S^2)$ and nontrivial harmonic maps $\{\omega_j\}_{j=1}^K$ for some integer $K>0$ such that $u(t_i)\rightarrow 0$ in $H^1(\S^2)$, $d(t_i)\rightharpoonup d_{\infty}$ in $H^1(\S^2,\S^2)$ and
\beq\label{engid}
\lim\limits_{t_i\rightarrow \infty}\int_{\S^2}|\nabla d(t_i)|^2\,dv_{0}
=\int_{\S^2}|\nabla d_{\infty}|^2\,dv_{0}+\sum\limits_{j=1}^K\int_{\S^2}|\nabla \omega_j|^2\,dv_{0}
\eeq

\item[(3)] If the initial data satisfies the following assumption
\beq\notag
\int_{\S^2}\left(|u_0|^2+|\nabla d_0|^2\right)\,dv_0(g)\leq 8\pi,
\eeq
then $u, d\in C^\infty\left(\S^2\times(0,\infty)\right)$. Moreover, there exist $t_i\uparrow \infty$ and a harmonic map $d_{\infty}\in C^\infty(\S^2,\S^2)$ such that
 $u(t_i)\rightarrow 0$ in $H^1(\S^2)$, $d(t_i)\rightarrow d_{\infty}$ in $H^1(\S^2,\S^2)$.
\end{itemize}
\end{theorem}

We would like to point out that the proof of Theorem \ref{thm1} should be slight modification of those in
\cite{huanglinwang14}, \cite{LTX16}, \cite{wangwang14}. The only difference is that we need to prove the energy inequality (see Lemma \ref{lemma0} below) with a weaker assumption \eqref{mus} for Leslie's coefficients instead of the following stronger one
\beq\label{nec2}
\mu_1+\frac{\lambda_2^2}{\lambda_1}\geq 0,\quad
\mu_4>0, \quad \mu_5+\mu_6\geq \frac{\lambda_2^2}{\lambda_1}.
\eeq
It is easy to see that the assumptions \eqref{parodi}, \eqref{nec1} and \eqref{nec2} imply \eqref{parodi}, \eqref{nec1} and \eqref{mus}.

Motivated  by the uniqueness results in\cite{wangxu14}  for the simplified system, we also show the uniqueness of the limit at time infinity to the initial value problem of the Ericksen-Leslie system \eqref{EL} and \eqref{initial}.
Before stating the main results, we need to recall some notations that have been introduced in Topping \cite{topping97}.

We define $z=x+iy$ as the complex coordinate on $\S^2$ which is homomorphic to $\bar{\mathbb C}$ via the stereographic projection. And we use the notation $dz=dx+idy$ and $d\z=dx-idy$. Then the metric $g_0$ on $\S^2$ can be written as $\sigma(z)^2dzd\z$ with
$$
\sigma(z)=\frac{2}{1+|z|^2},\quad z\in\bar{\mathbb C}.
$$
Similarly we assume that $v$ is the complex coordinate on the target $\S^2$ and the metric on the target $\S^2$ is $\rho(v)^2dvd\bar v$. For any $d\in H^1(\S^2,\S^2)$, we denote
$$
d_z=\frac12(d_x-id_y),\quad d_{\z}=\frac12(d_x+id_y).
$$
The $\partial$-energy and $\bar\partial$-energy of $d$ are respectively  given as follows
\beq\notag
E_{\partial}(d)=\frac{i}{2}\int_{\mathbb C}\rho^2(d)|d_z|^2\, dz\wedge d\z,\quad
E_{\bar\partial}(d)=\frac{i}{2}\int_{\mathbb C}\rho^2(d)|d_{\z}|^2\, dz\wedge d\z
\eeq
The Dirichlet energy of $d$ is given by
\beq\notag
E(d)=\frac12\int_{\S^2}|\nabla d|^2\,dv_0.
\eeq
It is not hard to see that
\beq\label{Ebar3}
E(d)=E_{\partial}(d)+E_{\bar\partial}(d),\quad
4\pi \mbox{deg}(d)=E_{\partial}(d)-E_{\bar\partial}(d),
\eeq
where $\mbox{deg}(d)$ is the topological degree of $d$, which is well-defined for $d\in H^1(\S^2,\S^2)$ (cf. \cite{BN95}).

Inspired by the conditions on the simplified system in \cite{wangxu14} (see also \cite{topping97} for similar results on heat flow of harmonic maps),  we are ready to state the sufficient conditions on the uniform limit at time infinity in $L^2$ under the weaker  parameter relationship.

\begin{theorem}\label{thm2}
If there exist $\ve_0>0$ and $T_0\geq1$ such that for the global weak solution $u\in L^{\infty}([0,\infty), L^2(\S^2,T\S^2))\cap L^{2}([0,\infty), H^1(\S^2,T\S^2))$, $d\in L^{\infty}([0,\infty), H^1(\S^2,\S^2))$ to the initial value problem \eqref{EL} and \eqref{initial} as in Theorem \ref{thm1}, it holds
\beq\label{assp1}
\frac{1}{2}\|u(T_0)\|_{L^2}^2+2\min \{E_{\partial} (d(T_0)),E_{\bar\partial} (d(T_0)) \}\leq \ve_0,
\eeq
then we should have as $t\rightarrow \infty$,
$u(t)\rightarrow 0$ in $H^1(\S^2)$, $d(t)\rightharpoonup d_{\infty}$ in $H^1(\S^2,\S^2)$ and $d(t)\rightarrow d_{\infty}$ in $L^2(\S^2,\S^2)$.
Furthermore, there exist integer $M_0>0$ and constants $C_1>0$, $C_2>0$ such that for any $t\geq T_0$, it holds
\beq\label{rate1}
\|u(t)\|_{L^2}+\|d(t)-d_{\infty}\|_{L^2}\leq C_1e^{-C_2t},
\eeq
\beq\label{rate2}
|E(d(t))-E(d_{\infty})-4M_0\pi|\leq C_1e^{-C_2t}.
\eeq
\end{theorem}

We also obtain the following result on the uniform limit at time infinity in $H^k$ with $k\geq 1$.

\begin{theorem}\label{thm3}
For the global weak solution $u\in L^{\infty}([0,\infty), L^2(\S^2,T\S^2))\cap L^{2}([0,\infty), H^1(\S^2,T\S^2))$, $d\in L^{\infty}([0,\infty), H^1(\S^2,\S^2))$ to the initial value problem \eqref{EL} and \eqref{initial} as in Theorem \ref{thm1}, suppose that there exists a time sequence $t_i\rightarrow \infty$ such that
 $u(t_i)\rightarrow 0$ in $L^2(\S^2)$, $d(t_i)\rightarrow d_{\infty}$ in $H^1(\S^2,\S^2)$ for some smooth harmonic map $d_{\infty}\in C^{\infty}(\S^2, \S^2)$. Then for any $k\geq1$, it holds $u(t)\rightarrow 0$ in $H^k(\S^2)$, $d(t)\rightarrow d_{\infty}$ in $H^k(\S^2,\S^2)$ as $t\rightarrow \infty$. Furthermore, there exist constants $C_1>0$ and $C_2>0$ depending only on $k$ such that
 \beq\label{rate3}
\|u(t)\|_{H^k}+\|d(t)-d_{\infty}\|_{H^k}\leq C_1e^{-C_2t}.
 \eeq
\end{theorem}

It is not hard to see from Theorem \ref{thm3} that the following assumption on initial data is also sufficient to obtain the uniform limit in $H^k$
 \beq\notag
 \frac{1}{2}\|u_0\|^2_{L^2}+E(d_0)\leq 4\pi,
 \eeq
 which is the assumption in part (3) of Theorem \ref{thm1}.


\section{Some Estimates}
\setcounter{equation}{0}


Firstly, under the weaker assumption \eqref{mus}, we give  the following apriori energy inequality for the regular solutions to the system \eqref{EL}. Similar arguments in
\cite{huanglinwang14}, \cite{LTX16}, \cite{wangwang14}, we should provide the proof of Theorem \ref{thm1}.


\begin{lemma}\label{lemma0}
Suppose that $(u,d)$ is a regular solution to the initial value problem \eqref{EL} and \eqref{initial} with \eqref{parodi}, \eqref{nec1} and \eqref{mus}. For any $t\in(0,\infty)$, the following energy inequality holds
\beq\label{enginq}
\frac{d}{dt}\frac{1}{2}\int_{\S^2}\left(|u|^2+|\nabla d|^2\right)\,dv_0
\leq-\int_{\S^2}\left(\alpha_0|\nabla u|^2+\frac{1}{\lambda_1}\left|\Delta d+|\nabla d|^2d\right|^2\right)\,dv_0,
\eeq
for some $\alpha_0>0$.
\end{lemma}

\pf
We multiplying the first equation in \eqref{EL} by $u$ and ues the integration by parts over $\S^2$ to obtain
\beq\label{inqpf1}
\begin{split}
\frac{d}{dt}\int_{\S^2}\frac12|u|^2\,dv_0
=-\int_{\S^2}(u\nabla d)\cdot\Delta d\,dv_0
-\int_{\S^2}\sigma^{L}:\nabla u\,dv_0.
\end{split}
\eeq
By the definition of $\sigma^L$, the second term on the right side can be computed as follows
\beq\label{inqpf2}
\begin{split}
&\sigma^{L}:\nabla u
=\sigma^{L}:(A+\Omega)\\
=&\mu_1\big(d^TAd\big)^2+\mu_4|A|^2+(\mu_5+\mu_6)|Ad|^2+\lambda_2N^TAd+\lambda_1d^T\Omega N+\lambda_2(d^T\Omega)(Ad).
\end{split}
\eeq
Putting \eqref{inqpf2} into \eqref{inqpf1} results in
\beq\label{inqpf3}
\begin{split}
&\frac{d}{dt}\int_{\S^2}\frac12|u|^2\,dv_0\\
=&-\int_{\S^2}(u\nabla d)\cdot\Delta d\,dv_0
-\int_{\S^2}\left(\mu_1\big(d^TAd\big)^2+\mu_4|A|^2+(\mu_5+\mu_6)|Ad|^2\right)\,dv_0\\
&-\int_{\S^2}\left(\lambda_2N^TAd+\lambda_1d^T\Omega N+\lambda_2(d^T\Omega)(Ad)\right)\,dv_0.
\end{split}
\eeq
Multiplying the third equation in \eqref{EL} by $N$ and integrating the resulting equation over $\S^2$, we can write
\beq\label{inqpf4}
\begin{split}
\frac{d}{dt}\int_{\S^2}\frac12|\nabla d|^2\,dv_0
=\int_{\S^2}(u\nabla d)\cdot\Delta d\,dv_0
-\int_{\S^2}\lambda_1N^2+\Delta d\cdot\Omega d+\lambda_2Ad\cdot N\,dv_0.
\end{split}
\eeq
Adding \eqref{inqpf3} to \eqref{inqpf4}, we are led to
\beq\label{inqpf5}
\begin{split}
&\frac{d}{dt}\int_{\S^2}\left(\frac12|u|^2+\frac12|\nabla d|^2\right)\,dv_0\\
=&
-\int_{\S^2}\left(\mu_1\big(d^TAd\big)^2+\mu_4|A|^2+(\mu_5+\mu_6)|Ad|^2+\lambda_1|N|^2\right)\,dv_0\\
&-\int_{\S^2}\left(2\lambda_2N^TAd+\lambda_1d^T\Omega N+\lambda_2(d^T\Omega)(Ad)+\Delta d\cdot \Omega d\right)\,dv_0.
\end{split}
\eeq
Since
\beq\notag
\lambda_1d^T\Omega N+\lambda_2(d^T\Omega)(Ad)
=(-\Omega d)^T\big(\lambda_1 N+\lambda_2(Ad)\big)
=(-\Omega d)\cdot \Delta d,
\eeq
we obtain
\beq\label{inqpf6}
\begin{split}
&\frac{d}{dt}\int_{\S^2}\left(\frac12|u|^2+\frac12|\nabla d|^2\right)\,dv_0\\
=&
-\int_{\S^2}\left(\mu_1\big(d^TAd\big)^2+\mu_4|A|^2+(\mu_5+\mu_6)|Ad|^2+\lambda_1|N|^2+2\lambda_2N^TAd\right)\,dv_0.
\end{split}
\eeq

Inspired by the arguments in \cite{AnnaLiu19} (see also \cite{WZZ13}), we can always find an orthonormal basis $\{e_1, e_2, e_3\}$ of $\R^3$ at any fixed $(x,t)$ such that
\beq
d=e_1,\quad N=\N e_2,\quad A=D_{ij}e_i\otimes e_j,
\eeq
where $\N:\S^2\times(0,\infty)\rightarrow \R$ and $D_{ij}:\S^2\times(0,\infty)\rightarrow\R$ satisfying $D_{ij}=D_{ji}$.
Then the right side of \eqref{inqpf6} can be written as
\beq\label{inqpf7}
\begin{split}
&\mu_1\big(d^TAd\big)^2+\mu_4|A|^2+(\mu_5+\mu_6)|Ad|^2+\lambda_1|N|^2+2\lambda_2N^TAd\\
=&\mu_1\big(D_{11}\big)^2+\mu_4\sum\limits_{ i, j=1}^3D_{ij}D_{ij}+(\mu_5+\mu_6)\sum\limits_{ i=1}^3|D_{1i}|^2+\lambda_1|\N|^2+2\lambda_2\N D_{12}.
\end{split}
\eeq

\smallskip
\noindent{\it Claim 1.\quad The dissipation terms in \eqref{inqpf7} are nonnegative with assumptions \eqref{parodi}, \eqref{nec1} and \eqref{mus}.}

\medskip
\noindent Indeed, the terms related to $\N$ and $D_{12}$ can be estimated as follows
\beq\label{inqpf8}
\begin{split}
&2\mu_4|D_{12}|^2+(\mu_5+\mu_6)|D_{12}|^2+\lambda_1|\N|^2+2\lambda_2\N D_{12}\\
=&(2\mu_4+\mu_5+\mu_6)|D_{12}|^2-\frac{\lambda_2^2}{\lambda_1}|D_{12}|^2+\left(\sqrt{\lambda_1}\N+\frac{\lambda_2}{\sqrt{\lambda_1}}D_{12}\right)^2\geq0,
\end{split}
\eeq
where we have used $\lambda_1>0$ and $(2\mu_4+\mu_5+\mu_6)>\frac{\lambda_2^2}{\lambda_1}$. Since $(2\mu_4+\mu_5+\mu_6)>\frac{\lambda_2^2}{\lambda_1}>0$ and $\mu_4>0$, it holds
\beq\label{inqpf9}
(2\mu_4+\mu_5+\mu_6)|D_{13}|^2\geq 0,\quad 2\mu_4|D_{23}|^2\geq 0.
\eeq

To control the rest terms in \eqref{inqpf7}, noting that
\beq\notag
\tr A=\nabla\cdot u=0,
\eeq
by which we obtain $\tr(D)=D_{11}+D_{22}+D_{33}=0$.
Then
\beq
|D_{33}|^2=|D_{11}|^2+|D_{22}|^2+2D_{11}D_{22}.
\eeq
Thus, the rest terms in \eqref{inqpf7} can be written as follows
\beq\label{inqpf10}
\begin{split}
&(\mu_1+\mu_4+\mu_5+\mu_6)|D_{11}|^2+\mu_4\left(|D_{22}|^2+|D_{33}|^2\right)\\
=&(\mu_1+2\mu_4+\mu_5+\mu_6)|D_{11}|^2+2\mu_4|D_{22}|^2+2\mu_4D_{11}D_{22}.
\end{split}
\eeq
The coefficient matrix of this quadratic form is given as follows
\beq
\left(
\begin{array}{cc}
\mu_1+2\mu_4+\mu_5+\mu_6&\mu_4\\
\\
\mu_4& 2\mu_4
\end{array}
\right).
\eeq
This matrix is positive definite if
\beq
\mu_1+2\mu_4+\mu_5+\mu_6>0
\eeq
and
\beq
2\mu_4(\mu_1+2\mu_4+\mu_5+\mu_6)-\mu_4^2>0.
\eeq
Combining the two inequalities above, we conclude that
\beq
2\mu_1+3\mu_4+2\mu_5+2\mu_6>0.
\eeq
Then \eqref{inqpf10} can be written as follows
\beq\label{inqpf11}
\begin{split}
&(\mu_1+\mu_4+\mu_5+\mu_6)|D_{11}|^2+\mu_4\left(|D_{22}|^2+|D_{33}|^2\right)\\
=&\left(\mu_1+\frac32\mu_4+\mu_5+\mu_6\right)|D_{11}|^2+\mu_4\left(\frac{1}{\sqrt{2}}D_{11}+\sqrt{2}D_{22}\right)^2\geq 0,
\end{split}
\eeq
which combining with \eqref{inqpf8} and \eqref{inqpf9} indicates Claim 1.

\medskip
\noindent{\it Claim 2. \quad There exists a positive constant $\alpha_0>0$ such that it holds
 \beq\label{inqpf0}
\begin{split}
&\mu_1\big(D_{11}\big)^2+\mu_4\sum\limits_{ i, j=1}^3D_{ij}D_{ij}+(\mu_5+\mu_6)\sum\limits_{ i=1}^3|D_{1i}|^2+\lambda_1|\N|^2+2\lambda_2\N D_{12}\\
\geq&\frac{1}{\lambda_1}|\Delta d+|\nabla d|^2d|^2+2\alpha_0|A|^2.
\end{split}
\eeq
}
\smallskip

\noindent Indeed, by the proof of Claim 1, we obtain
\beq\label{inqpf15}
\begin{split}
&\mu_1\big(D_{11}\big)^2+\mu_4\sum\limits_{ i, j=1}^3D_{ij}D_{ij}+(\mu_5+\mu_6)\sum\limits_{ i=1}^3|D_{1i}|^2+\lambda_1|\N|^2+2\lambda_2\N D_{12}\\
=&(2\mu_4+\mu_5+\mu_6)|D_{12}|^2-\frac{\lambda_2^2}{\lambda_1}|D_{12}|^2+\left(\sqrt{\lambda_1}\N+\frac{\lambda_2}{\sqrt{\lambda_1}}D_{12}\right)^2\\
&+(\mu_1+\mu_4+\mu_5+\mu_6)|D_{11}|^2+(2\mu_4+\mu_5+\mu_6)|D_{13}|^2+\mu_4\left(|D_{22}|^2+2|D_{23}|^2+|D_{33}|^2\right).
\end{split}
\eeq
Direct computation implies
\beq\notag
\begin{split}
&|\Delta d+|\nabla d|^2d|^2=\left|\lambda_1 N+\lambda_2Ad\right|^2-\lambda_2^2|d^TAd|^2\\
=&\lambda_1^2\N^2 +\lambda_2^2\sum\limits_{ i=1}^3|D_{1i}|^2+2\lambda_1\lambda_2\N D_{12}-\lambda_2^2|D_{11}|^2\\
=&\left(\lambda_1\N+\lambda_2D_{12}\right)^2+\lambda_2^2|D_{13}|^2
\end{split}
\eeq
from which we have
\beq\notag
\begin{split}
&(2\mu_4+\mu_5+\mu_6)|D_{13}|^2+\left(\sqrt{\lambda_1}\N+\frac{\lambda_2}{\sqrt{\lambda_1}}D_{12}\right)^2\\
=&(2\mu_4+\mu_5+\mu_6)|D_{13}|^2+\frac{1}{\lambda_1}\left(\lambda_1\N+\lambda_2D_{12}\right)^2\\
=&\left(2\mu_4+\mu_5+\mu_6-\frac{\lambda_2^2}{\lambda_1}\right)|D_{13}|^2+\frac{1}{\lambda_1}|\Delta d+|\nabla d|^2d|^2.
\end{split}
\eeq
It is easy to see that
\beq
|A|^2=\sum\limits_{i,j=1}^3|D_{ij}|^2.
\eeq
Choose
\beq
\alpha_0=\frac14\min\left\{2\mu_4+\mu_5+\mu_6-\frac{\lambda_2^2}{\lambda_1},\ 2\mu_4,\ \frac{\mu_1+4\mu_4+\mu_5+\mu_6-\Delta_0}{8}\right\}
\eeq
where
\beq
\Delta_0=\sqrt{(\mu_1+\mu_5+\mu_6)^2+4\mu_4^2}.
\eeq
By the assumption \eqref{mus}, direct computation implies
\beq
\gamma_1:=\frac{\mu_1+4\mu_4+\mu_5+\mu_6-\Delta_0}{2}>0.
\eeq
Here $\gamma_1$ is the smaller eigenvalue of the coefficient matrix of the quadratic form in \eqref{inqpf10}.
Thus the dissipation terms in \eqref{inqpf15} can be estimates as follows
\beq\label{inqpf12}
\begin{split}
&\mu_1\big(D_{11}\big)^2+\mu_4\sum\limits_{ i, j=1}^3D_{ij}D_{ij}+(\mu_5+\mu_6)\sum\limits_{ i=1}^3|D_{1i}|^2+\lambda_1|\N|^2+2\lambda_2\N D_{12}\\
=&\left(2\mu_4+\mu_5+\mu_6-\frac{\lambda_2^2}{\lambda_1}\right)\left(|D_{12}|^2+|D_{13}|^2\right)+\frac{1}{\lambda_1}|\Delta d+|\nabla d|^2d|^2+2\mu_4|D_{23}|^2\\
&+(\mu_1+\mu_4+\mu_5+\mu_6)|D_{11}|^2+\mu_4\left(|D_{22}|^2+|D_{33}|^2\right)\\
\geq&\frac{1}{\lambda_1}|\Delta d+|\nabla d|^2d|^2+4\alpha_0\left(|D_{12}|^2+|D_{13}|^2+|D_{23}|^2\right)+\gamma_1\left(|D_{11}|^2+|D_{22}|^2\right).
\end{split}
\eeq
For the last term, by the traceless condition of $D$, it holds
\beq
|D_{11}|^2+|D_{22}|^2\geq \frac12\left(|D_{11}|^2+|D_{22}|^2\right)+\frac14|D_{33}|^2
\eeq
Putting this into \eqref{inqpf12}, we conclude that
\beq\label{inqpf13}
\begin{split}
&\mu_1\big(D_{11}\big)^2+\mu_4\sum\limits_{ i, j=1}^3D_{ij}D_{ij}+(\mu_5+\mu_6)\sum\limits_{ i=1}^3|D_{1i}|^2+\lambda_1|\N|^2+2\lambda_2\N D_{12}\\
\geq&\frac{1}{\lambda_1}|\Delta d+|\nabla d|^2d|^2+4\alpha_0\left(|D_{12}|^2+|D_{13}|^2+|D_{23}|^2+|D_{11}|^2+|D_{22}|^2+|D_{33}|^2\right)\\
\geq&\frac{1}{\lambda_1}|\Delta d+|\nabla d|^2d|^2+2\alpha_0|A|^2,
\end{split}
\eeq
which indicates Claim 2.

\medskip
Plugging \eqref{inqpf13} into \eqref{inqpf6}, we conclude
\beq\label{inqpf14}
\begin{split}
&\frac{d}{dt}\int_{\S^2}\left(\frac12|u|^2+\frac12|\nabla d|^2\right)\,dv_0\\
\leq&
-\int_{\S^2}\left(\frac{1}{\lambda_1}|\Delta d+|\nabla d|^2d|^2+2\alpha_0|A|^2\right)\,dv_0\\
=&-\int_{\S^2}\left(\frac{1}{\lambda_1}|\Delta d+|\nabla d|^2d|^2+\alpha_0|\nabla u|^2\right)\,dv_0
\end{split}
\eeq
where we have used the fact $\nabla\cdot u=0$ and hence
\beq
\int_{\S^2}|A|^2\,dv_0=\frac12\int_{\S^2}|\nabla u|^2\,dv_0.
\eeq

\endpf



Similar to the proof of Lemma \ref{lemma0}, we are able to show the following energy estimates for the global weak solution obtained in Theorem \ref{thm1}.

\begin{lemma}\label{lemma1}
Suppose that $u\in L^{\infty}([0,\infty), L^2(\S^2,T\S^2))\cap L^{2}([0,\infty), H^1(\S^2,T\S^2))$,  $d\in L^{\infty}([0,\infty), \\H^1(\S^2,\S^2))$ is the global weak solution to the initial value problem \eqref{EL} and \eqref{initial} obtained in Theorem \ref{thm1}. For any $t\in (0,\infty)\setminus\{T_l\}_{l=1}^L$, the following energy equality holds
\beq\label{engest}
\frac{d}{dt}\frac{1}{2}\int_{\S^2}\left(|u|^2+|\nabla d|^2\right)\,dv_0
\leq-\int_{\S^2}\left(\alpha_0|\nabla u|^2+\frac{1}{\lambda_1}\left|\Delta d+|\nabla d|^2d\right|^2\right)\,dv_0.
\eeq
\end{lemma}

We also need the following estimate for the global weak solution obtained in Theorem \ref{thm1}.
\begin{lemma}\label{lemma2}
Suppose that $u\in L^{\infty}([0,\infty), L^2(\S^2,T\S^2))\cap L^{2}([0,\infty), H^1(\S^2,T\S^2))$, $d\in L^{\infty}([0,\infty), \\H^1(\S^2,\S^2))$ is the global weak solution to the initial value problem \eqref{EL} and \eqref{initial} obtained in Theorem \ref{thm1}. For any $t\in [0,\infty)\setminus\{T_l\}_{l=1}^L$, it holds
\beq\label{uavg0}
\int_{\S^2}u(x,t)\,dv_0=0.
\eeq
\end{lemma}

\pf
This can be proved by integrating the first equation of \eqref{EL} over $\S^2$ and using the assumption \eqref{u0con} on initial data.

\endpf

Finally, we state the key estimates of $E_{\partial}$ or $E_{\bar\partial}$ in terms of tension field, which have been first proved in Topping \cite{topping97}.
\begin{lemma}\label{lemma3}
There exist constants $\ve_0>0$ and $C>0$ such that for any $d\in H^1(\S^2,\S^2)$
\begin{itemize}
\item[(1)] If it holds $E_{\partial}(d)<\ve_0$,
then
\beq\label{Ebar2}
E_{\partial}(d)
\leq C\int_{\S^2}\big|\Delta d+|\nabla d|^2\big|^2\,dv_0.
\eeq

\item[(2)] If it holds
$E_{\bar\partial}(d)<\ve_0$,
then
\beq\label{Ebar1}
E_{\bar\partial}(d)
\leq C\int_{\S^2}\big|\Delta d+|\nabla d|^2\big|^2\,dv_0.
\eeq
\end{itemize}
\end{lemma}

%


\section{Long Time Dynamics}
\setcounter{equation}{0}


 In this section, we will devote to the proof of the uniform limit at time infinity in $L^2$ and $H^k$ with the help of the energy estimates.

 \vspace{3mm}
\noindent{\bf Proof of Theorem \ref{thm2}}.\quad
Firstly, in order to proof the strong convergence $u(t)\rightarrow 0$ in $H^1(\S^2)$, it is sufficient to show that $\|u\|_{H^{1}}$ is exponentially decaying.
For any $t>T_0$, we integrate from \eqref{engest} over $(T_0, t)$ in Lemma \ref{lemma1} to obtain
\beq\label{engt0}
\frac{1}{2}\int_{\S^2}\left(|u(t)|^2+|\nabla d(t)|^2\right)\,dv_0
\leq \frac{1}{2}\int_{\S^2}\left(|u(T_0)|^2+|\nabla d(T_0)|^2\right)\,dv_0.
\eeq
Taking $\ve_0$ as in Lemma \ref{lemma3} and by the assumption \eqref{assp1}, we should have
\beq\label{assp2}
\min\big\{E_{\partial}(d(T_0)),E_{\bar\partial}(d(T_0))\big\}<\ve_0.
\eeq
Since $\mbox{deg} (d)$ is a constant for any $t\geq T_0$, by the energy estimate \eqref{engt0} and the relation \eqref{Ebar3}, we deduce that
\beq\label{engt1}
\frac{1}{2}\|u(t)\|^2_{L^2}+2E_{\partial}(d(t))
\leq \frac{1}{2}\|u(T_0)\|^2_{L^2}+2E_{\partial}(d(T_0)),
\eeq
together with assumption \eqref{assp1}, which is a uniform upper bound
for any $t\geq T_0$.
By Lemma \ref{lemma2},  Lemma \ref{lemma3} and the Poincar\'e inequality, we have
\beq\label{engt4}
\frac{1}{2}\|u(t)\|_{L^2}^2+2E_{\partial} (d(t))
\leq C\|\nabla u(t)\|_{L^2}^2+C\|\Delta d(t)+|\nabla d(t)|^2d(t)\|_{L^2}^2.
\eeq
By the identity \eqref{Ebar3} and \eqref{engest}, it is not hard to see that for any $t\geq T_0$
\beq\label{engt5}
\begin{split}
\frac{d}{dt}\left(\frac{1}{2}\|u\|_{L^2}^2+2E_{\partial} (d)\right)
\leq-\int_{\S^2}\left(\alpha_0|\nabla u|^2+\frac{1}{\lambda_1}\left|\Delta d+|\nabla d|^2d\right|^2\right)\,dv_0,
\end{split}
\eeq
Combining \eqref{engt4} with \eqref{engt5} yields
\beq\label{engt7}
\begin{split}
-\frac{d}{dt}\left(\frac{1}{2}\|u\|_{L^2}^2+2E_{\partial} (d)\right)^{\frac12}
\geq &C \left(\|\nabla u\|_{L^2}^2+\|\Delta d+|\nabla d|^2d\|_{L^2}^2\right)^{\frac12}\\
\geq &C \left(\frac{1}{2}\|u\|_{L^2}^2+2E_{\partial} (d)\right)^{\frac12}
\end{split}
\eeq
where the constant $C$ depends on $\min\{\alpha_0, \frac{1}{\lambda_1}\}$. We apply the Gronwall inequality to deduce
\beq\label{engt8}
\|u(t)\|_{L^2}^2+E_{\partial} (d(t))
\leq \left(\|u(T_0)\|_{L^2}^2+E_{\partial} (d(T_0))\right)e^{-C(t-T_0)}
\leq \ve_0 e^{-C(t-T_0)}.
\eeq
Next, we consider  $d(t)\rightarrow d_{\infty}$ in the sense of $L^2(\S^2,\S^2)$. We integrate the first inequality in \eqref{engt7} for any $t\geq 2T_0$ to obtain
\beq\label{engt9}
\int_{t}^\infty\left(\|\nabla u\|_{L^2}^2+\|\Delta d+|\nabla d|^2d\|_{L^2}^2\right)^{\frac12}\,ds
\leq \left(\frac{1}{2}\|u(t)\|_{L^2}^2+2E_{\partial} (d(t))\right)^{\frac12}\leq Ce^{-Ct}.
\eeq
By the equation of $d$ in \eqref{EL}, we obtain for any $\tilde t\geq t\geq 2T_0$
\beq\label{engt10}
\begin{split}
&\|d(\tilde t)-d(t)\|_{L^2}^2
\leq 2\|d(\tilde t)-d(t)\|_{L^1}
\leq 2\int_{t}^{\tilde t}\|d_{t}\|_{L^1}\,ds\\
\leq& C\int_{t}^{\tilde t}\|u\|_{L^2}\|\nabla d\|_{L^2}\,ds+C\int_{t}^{\tilde t}\|\nabla u\|_{L^1}\,ds+C\int_{t}^{\tilde t}\|\Delta d+|\nabla d|^2d\|_{L^1}\,ds\\
\leq &C\int_{t}^{\tilde t}\|\nabla u\|_{L^2}\,ds+C\int_{t}^{\tilde t}\|\Delta d+|\nabla d|^2d\|_{L^2}\,ds\\
\leq & C\left(\|u(t)\|_{L^2}^2+E_{\partial} (d(t))\right)^{\frac12}\leq Ce^{-Ct},
\end{split}
\eeq
where we have also used \eqref{engt0} and \eqref{engt9} in the forth and fifth inequalities respectively. It is not hard to see that the estimate \eqref{engt10} implies that $d(t)$ is convergent in $L^2$ as $t\rightarrow \infty$. By the conclusion in Theorem \ref{thm1}, there exists a smooth harmonic map $d_\infty\in C^{\infty}(\S^2, \S^2)$ and a sequence $t_i $ such that $d(t_i)\rightharpoonup d_{\infty}$ in $H^1(\S^2,\S^2)$. Taking $\tilde t=t_i$ and letting $i\rightarrow\infty$, it holds for any $t\geq 2T_0$
\beq\label{engt11}
\|d(t)-d_{\infty}\|_{L^2}^2\leq Ce^{-Ct}.
\eeq
Let $\tilde t_i\geq 2T_0$ be any sequence that tends to $\infty$. By the energy estimate \eqref{engt0}, there exists a subsequence $\tilde t_{i_j}\rightarrow \infty$ such that  $d(\tilde t_{i_j})\rightharpoonup \tilde d_{\infty}$ in $H^1(\S^2,\S^2)$. By \eqref{engt11}, we should have $d_{\infty}=\tilde d_{\infty}$, which implies $d(t)\rightharpoonup d_{\infty}$ in $H^1(\S^2,\S^2)$ as $t\rightarrow \infty$.

To this end, we  calculate for any $\tilde t\geq t\geq 2T_0$
\beq\notag
\begin{split}
0\leq& \left(\frac{1}{2}\|u(t)\|^2_{L^2}+E(d(t))\right)
-\left( \frac{1}{2}\|u(\tilde t)\|^2_{L^2}+E(d(\tilde t))\right)\\
\leq &\left(\frac{1}{2}\|u(t)\|^2_{L^2}+2E_{\partial}(d(t))\right)
-\left( \frac{1}{2}\|u(\tilde t)\|^2_{L^2}+2E_{\partial}(d(\tilde t))\right)\\
\leq &\left(\frac{1}{2}\|u(t)\|^2_{L^2}+2E_{\partial}(d(t))\right)\leq Ce^{-Ct},
\end{split}
\eeq
which is similar to the derivation of \eqref{engt0} and \eqref{engt1}. Thus we lead to
\beq\label{engt12}
|E(d(t))-E(d(\tilde t))|\leq \frac{1}{2}|\|u(t)\|_{L^2}^2-\|u(\tilde t)\|_{L^2}^2|+Ce^{-Ct}\leq Ce^{-Ct},
\eeq
where we have used the fact \eqref{engt8} in last inequality.  It is not hard to see that the estimate \eqref{engt12} implies that $E(d(t))$ converges as $t\rightarrow \infty$. By the conclusion in Theorem \ref{thm1}, there exists a smooth harmonic map $d_\infty\in C^{\infty}(\S^2, \S^2)$, nontrivial harmonic maps $\{\omega_j\}_{j=1}^K$and a sequence $t_i\rightarrow \infty$ such that
\eqref{engid} holds. Since $\omega_j$ is nontrivial harmonic maps, there are positive integer $m_j$ such that
$$
\frac{1}{2}\int_{\S^2}|\nabla \omega_j|^2\,dv_0=4m_j\pi.
$$
Choosing $\tilde t=t_i$ and letting $i\rightarrow \infty$, we obtain
$$
E(d(\tilde t))\rightarrow E(d_{\infty})+4M_0\pi
$$
for some nonnegative integer $M_0$. Putting the preceding limit into \eqref{engt12} yields
\beq
|E(d(t))-E(d_{\infty})-4M_0\pi|\leq  Ce^{-Ct}
\eeq
for any $t\geq 2T_0$. Thus, we obtain the desired result.

\endpf

By Theorem \ref{thm2}, we have the following result.
\begin{corollary}\label{thm4}
Suppose that $u\in L^{\infty}([0,\infty), L^2(\S^2,T\S^2))\cap L^{2}([0,\infty), H^1(\S^2,T\S^2))$, $d\in L^{\infty}([0,\infty),\\ H^1(\S^2,\S^2))$
is the global weak solution to the initial value problem \eqref{EL} and \eqref{initial} as in Theorem \ref{thm1}.

\begin{itemize}

\item[(1)] Suppose that for a sequence $t_i\rightarrow \infty$, the harmonic maps $d_{\infty}$ and $\{\omega_j\}_{j=1}^K$ are the weak limit and bubbles associated to $t_i$, which are all holomorphic or all anti-holomorphic. Then as $t\rightarrow \infty$,
$u(t)\rightarrow 0$ in $H^1(\S^2)$, $d(t)\rightharpoonup d_{\infty}$ in $H^1(\S^2,\S^2)$ and $d(t)\rightarrow d_{\infty}$ in $L^2(\S^2,\S^2)$.
Furthermore, there exist integer $M_0>0$ and constants $C_1>0$, $C_2>0$ such that for any $t\geq T_0$, it holds
\beq\label{rate4}
\|u(t)\|_{L^2}+\|d(t)-d_{\infty}\|_{L^2}\leq C_1e^{-C_2t},
\eeq
\beq\label{rate5}
|E(d(t))-E(d_{\infty})-4M_0\pi|\leq C_1e^{-C_2t}.
\eeq

\item[(2)] If the initial data satisfies
\beq\label{initialsm}
\frac{1}{2}\|u_0\|_{L^2}^2+2\min \{E_{\partial} (d_0),E_{\bar\partial} (d_0) \}< 8\pi,
\eeq
then we obtain the same conclusions.
\end{itemize}
\end{corollary}

\pf For the proof of Part (1), we only need to confirm whether the assumption  \eqref{assp1} of Theorem \ref{thm2} is satisfied. Without loss of generality, we assume that $d_\infty$ and all $\{\omega_j\}_{j=1}^K$ are all anti-holomorphic. Then it holds for any $j=1,\ldots,K$
\beq\notag
E_{\partial}(d_\infty)=E_{\partial}(\omega_j)=0.
\eeq
By the limit \eqref{engid} in Theorem \ref{thm1}, we obtain as $t_i\rightarrow \infty$
\beq\notag
E_{\partial}(d(t_i))\rightarrow 0.
\eeq
Combining the fact $u(t_i)\rightarrow 0$ in $H^1(\S^2)$, we should be able to find a time $t_{0}$ large enough such that the assumption \eqref{assp1} holds.

In order to complete the proof of Part (2), we only need to verify the assumption of \eqref{initialsm} satisfying that the weak limit $d_{\infty}$ and all  bubbles $\{\omega\}_{j=1}^K$ are all holomorphic or all anti-holomorphic. Without loss of generality, we assume that
$$E_{\partial}(d_0)\leq E_{\bar\partial}(d_0).$$

Let $T_1>0$ be the first possible singular time of the system \eqref{EL}. we claim that
$\mbox{deg}(d(t))$ is constant for any $0\leq t<T_1$. In fact, it follows that the solution $(u,d)$ is continuous  before $T_1$ and let $t\rightarrow 0^+$,
$$
\|u(t)-u_0\|_{L^2}+\|\nabla d(t)-\nabla d_0\|_{L^2}\rightarrow 0.
$$
Suppose $\{T_l\}_{l=1}^L$ is all the possible singular times, making use of the fact \eqref{Ebar3} and integrating the energy law \eqref{engest} in Lemma \ref{lemma1} result in
\beq\notag
\frac{1}{2}\|u(t)\|_{L^2}^2+2E_{\partial} (d(t))\leq \frac{1}{2}\|u(T_l)\|_{L^2}^2+2E_{\partial} (d(T_l))
\eeq
for any $t\in [T_l, T_{l+1}]$ with $l=1, \ldots, L$ and $T_{L+1}:=\infty$ and $[T_L,T_{L+1}]:=[T_L,\infty)$.
Therefore for any $t\geq 0$, it holds
\beq\notag
\frac{1}{2}\|u(t)\|_{L^2}^2+2E_{\partial} (d(t))\leq \frac{1}{2}\|u_0\|_{L^2}^2+2E_{\partial} (d_0)<8\pi,
\eeq
where we have used the assumption \eqref{initialsm} in the last inequality. The lower semicontinuity implies that
\beq\notag
2E_{\partial}(d_\infty)+\sum_{l=1}^L2E_{\partial}(\omega_l)\leq \lim\limits_{t_i\rightarrow \infty} \left(\frac{1}{2}\|u(t)\|_{L^2}^2+2E_{\partial} (d(t))\right)<8\pi.
\eeq
We complete the proof.

\endpf

%
%

We devote the rest of this section to prove the uniform limit at time infinity in $H^k$.

 \vspace{3mm}
\noindent{\bf Proof of Theorem \ref{thm3}}.\quad

By the interpolation inequality, we claim the decay rate estimate \eqref{rate3}, provided the following inequalities
\beq\label{pfth31}
\|u(t)\|_{L^2}+\|d(t)-d_{\infty}\|_{L^2}\leq C_1e^{-C_2t}
\eeq and \beq\label{pfth34}
\sup\limits_{t\geq \tilde T}\left(\|u(t)\|_{C^k}+\| d(t)\|_{C^{k+1}}\right)\leq C_k
\eeq
hold. Without loss of generality, we may assume that $d_{\infty}$ is anti-holomorphic since all the harmonic maps from $\S^2$ to $\S^2$ are either holomorphic or anti-holomorphic. Direct computation implies
\beq\notag
E_{\partial}(d(t_i))=E_{\partial}(d(t_i)-d_\infty)\leq E(d(t_i)-d_{\infty}).
\eeq
The strong convergence of $u(t_i)$ and $\nabla d(t_i)$ implies that there exists a $T_0>1$ such that
\beq\notag
\frac12\|u(T_0)\|_{L^2}+2E_{\partial}(d(T_0))<\ve_0
\eeq
where $\ve_0$ is given in Lemma \ref{lemma3}. Theorem \ref{thm2} can be applied directly and it holds as $t\rightarrow \infty$ that $u(t)\rightarrow 0$ in $H^1(\S^2)$, $d(t)\rightharpoonup d_{\infty}$ in $H^1(\S^2,\S^2)$ and
\eqref{pfth31} is deduced.

The strong convergence of $u(t_i)$ and $\nabla d(t_i)$ implies
\beq\notag
\frac12\|u(t_i)\|_{L^2}+E(d(t_i))\rightarrow E(d_{\infty}).
\eeq
Combining with the energy decay estimate in Lemma \ref{lemma1}, it follows that
\beq\label{pfth32}
E(d(t))\rightarrow E(d_{\infty})
\eeq
as $t\rightarrow\infty$. For any $\ve>0$, we can find a constant $r_0>0$ such that
\beq\notag
\max\limits_{x\in \S^2}\int_{\S^2\cap B_{r_0}(x)}|\nabla d_{\infty}|^2\,dv_0\leq \frac{\ve^2}{2}.
\eeq
By the convergence \eqref{pfth32}, we can find a large time $\tilde T$ such that
\beq\label{pfth33}
\sup\limits_{t\geq \tilde T}\,\max\limits_{x\in \S^2}\int_{\S^2\cap B_{r_0}(x)}\left(|u(t)|^2+|\nabla d(t)|^2\right)\,dv_0\leq \ve^2.
\eeq
By the proof of Lemma 4.2 in \cite{huanglinwang14}, we can show the following estimate
\beq
\int_{\tilde T}^{\infty}\int_{\S^2}\left(|\nabla u|^2+|\nabla^2 d|^2\right)\,dv_0dt
\leq \tilde C,
\eeq
where the positive constant $\tilde C$ only depends on $\|u(\tilde T)\|_{L^2}$ and $E(d(\tilde T))$. By Theorem 1.3 in \cite{huanglinwang14}, we conclude the regularity of the solution. Furthermore, for any $k\geq 0$, we can conclude \eqref{pfth34}.
The proof is completed.

\endpf
\section*{Conflicts of Interest}
Authors have no conflict of interest to declare.

\section*{Funding}
This work was supported by the National Natural Science Foundation of China ( 12071058, 11871134).

\end{document}